\documentclass[10pt]{article}
\usepackage{amsmath,amsthm,amssymb,amsfonts,mathrsfs}
\usepackage{graphicx}
\usepackage[square,numbers,sort&compress]{natbib}

\ifx\pdfoutput\undefined
 \usepackage[dvipdfm,
  pdfstartview=FitH,
  bookmarks=true,
  bookmarksnumbered=true,
  bookmarksopen=true,
  plainpages=false,
  pdfpagelabels,
  colorlinks=true,
  linkcolor=black,
  citecolor=black,
  pdfborder=001]{hyperref}
\else
 \usepackage[pdftex,
  pdfstartview=FitH,
  bookmarks=true,
  bookmarksnumbered=true,
  bookmarksopen=true,
  plainpages=false,
  pdfpagelabels,
  colorlinks=true,
  linkcolor=black,
  citecolor=black,
  pdfborder=001]{hyperref}
\fi

\newtheorem{theorem}{Theorem}[section]
\newtheorem{lemma}{Lemma}

\newtheorem*{conjecture*}{Conjecture}
\newtheorem{claim}{Claim}
\theoremstyle{definition}
\newtheorem{remark}{Remark}

\begin{document}
\title{Matching and Factor-Critical Property in $3$-Dominating-Critical Graphs}

\date{}

\maketitle

\begin{abstract}
  Let $\gamma(G)$ be the domination number of a graph $G$.
A graph $G$ is \emph{domination-vertex-critical}, or
\emph{$\gamma$-vertex-critical}, if $\gamma(G-v)< \gamma(G)$ for
every vertex $v \in V(G)$.  In this paper, we show that:  Let $G$ be
a $\gamma$-vertex-critical graph and $\gamma(G)=3$. (1) If $G$ is of
even order and $K_{1,6}$-free, then $G$ has a perfect matching; (2)
If $G$ is of odd order and $K_{1,7}$-free, then  $G$ has a near
perfect matching with only three exceptions. All these results
improve the known results.\\


\end{abstract}

\section{Introduction}
Let $G$ be a finite simple graph with vertex set $V(G)$ and edge set
$E(G)$. A set $S\subseteq  V$ is a \emph{dominating set} of $G$ if
every vertex in $V$ is either in $S$ or is adjacent to a vertex in
$S$. For two sets $A$ and $B$, $A$ \emph{dominates} $B$ if every
vertex of $B$ has a neighbor in $A$ or is a vertex of $A$;
sometimes, we also say that $B$ is \emph{dominated} by $A$. Let
$u\in V$ and $A\subseteq  V-\{u\}$, if $u$ is adjacent to some
vertex of $A$, then we say that $u$ is \emph{adjacent} to $A$. The
\emph{domination number} of $G$, denoted by $\gamma(G)$, is the
minimum cardinality of dominating sets of $G$. A graph $G$ is
\emph{domination vertex critical}, or
\emph{$\gamma$-vertex-critical}, if $\gamma(G-v)< \gamma(G)$ for
every vertex $v \in V(G)$. Indeed, if $\gamma(G-v)<\gamma(G)$, then
$\gamma(G-v)=\gamma(G)-1$. A graph $G$ is \emph{domination edge
critical}, if $\gamma(G+e)<\gamma(G)$ for any edge $e\notin E(G)$.
We call a graph $G$ \emph{$k$-$\gamma$-vertex-critical} (resp.
\emph{$k$-$\gamma$-edge-critical}) if it is domination vertex
critical (resp. domination edge critical) and $\gamma(G)=k$.

    A matching is \emph{perfect} if it is incident with every vertex of
$G$. If $G-v$ has a perfect matching for every choice of $v\in
V(G)$, $G$ is said to be \emph{factor-critical}.  The concept of
factor-critical graphs was first introduced by Gallai in 1963 and it
plays an important role in the study of matching theory. Contrary to
its apparent strict condition, such graphs form a relatively rich
family for study. It is the essential ``building block" for
well-known Gallai-Edmonds Matching Structure Theorem.

    The subject of $\gamma$-vertex-critical graphs was studied first by
Brigham, Chinn and Dutton \cite{Brigham1988} and continued by Fulman,
Hanson and MacGillivray \cite{Fulman1995}. Clearly, the only
$1$-$\gamma$-vertex-critical graph is $K_1$ (i.e., a single vertex).
Brigham, Chinn and Dutton \cite {Brigham1988} pointed out that the
$2$-$\gamma$-vertex-critical graphs are precisely the family of
graphs obtained from the complete graphs $K_{2n}$ with a perfect
matching removed (Theorem \ref{2vertex-critical}). For $k>2$, however,
much remains unknown about the structure of
$k$-$\gamma$-vertex-critical graphs. Recently, Ananchuen and Plummer
\cite {Ananchuen2005, Ananchuen2007} began to investigate matchings
in $3$-$\gamma$-vertex-critical graphs. They showed that a
$K_{1,5}$-free $3$-$\gamma$-vertex-critical graph of {\it even} order
has a perfect matching (see \cite {Ananchuen2005}). For the graphs of
{\it odd} order, they proved that the condition of $K_{1,4}$-freedom
is sufficient for factor-criticality (see \cite {Ananchuen2007}).
Wang and Yu \cite{Wang2008} improved this result by weakening the
condition of $K_{1,4}$-freedom to almost $K_{1,5}$-freedom.  In
\cite{Wang}, they also studied the $k$-factor-criticality in
$3$-$\gamma$-edge-critical graphs and obtained several useful results
on connectivity of $3$-$\gamma$-vertex-critical graphs.

  The relevant theorems are stated formally below.

\begin{theorem}[Brigham {\it et al}., \cite{Brigham1988}]\label{2vertex-critical}
A graph $G$ is $2$-$\gamma$-vertex-critical if and only if it is
isomorphic to $K_{2n}$ with a perfect matching removed.
\end{theorem}

\begin{theorem}[Ananchuen and Plummer, \cite{Ananchuen2005}]\label{AP05}
  Let $G$ be a $3$-$\gamma$-vertex-critical graph of even order. If
  $G$ is $K_{1,5}$-free, then $G$ has a perfect matching.
\end{theorem}

\begin{theorem}[Ananchuen and Plummer, \cite{Ananchuen2007}]\label{AP07}
  Let $G$ be a $3$-$\gamma$-vertex-critical graph of odd order at
  least $11$. If $G$ is $K_{1,5}$-free, then $G$ contains a near
  perfect matching.
\end{theorem}

   For $v\in V(G)$, we denote a minimum dominating set of $G-v$ by
$D_v$. The following facts about $D_v$ follow immediately from the
definition of $3$-$\gamma$-vertex-criticality and we shall use it
frequently in the proofs of the main theorems.

\vspace{4mm}

 \noindent\textbf{Facts:} If $G$ is
$3$-$\gamma$-vertex-critical, then the followings hold
\begin{itemize}
\vspace{-5pt}
\item [(1)] For every vertex $v$ of $G$, $|D_v|=2$;
\vspace{-5pt}
\item [(2)] If $D_v=\{x,y\}$, then $x$ and $y$ are not adjacent to
$v$; \vspace{-5pt}
\item [(3)] For every pair of distinct vertices $v$ and $w$, $D_v\neq D_w$.
\end{itemize}

  In this paper, we utilize the techniques developed in \cite{Wang2008}
and \cite{Wang} to extend Theorem \ref{AP05} and Theorem \ref{AP07} to the
following theorem.

\begin{theorem}\label{k167_theorem}
  Let $G$ be a $3$-$\gamma$-vertex-critical graph.
\begin{enumerate}
  \item  If $G$ is $K_{1,6}$-free and $|V(G)|$ is even, $|V(G)|\neq 12$, then $G$ has a perfect
  matching.
  \item  If $G$ is $K_{1,7}$-free of odd order, and $c_o(G)=1$, $|V(G)| \neq 13$, then either $G$ has a
  near perfect matching or $G$ is one of Fig. \ref{9vertex} and Fig. \ref{s45}.
\end{enumerate}
\end{theorem}

In theory of matching, Tutte's $1$-Factor Theorem plays a central
role. From $1$-Factor Theorem, a characterization of a graph with a
near perfect matching can be easily derived. Following the
convention of \cite{Lov'asz1986}, we use $c(G)$ (resp. $c_o(G)$) to
denote the number of (resp. odd) components of $G$.

\begin{theorem}[Tutte's $1$-Factor Theorem]\label{tutte_thm}
  A graph $G$ has a perfect matching if and only if for any $S\subseteq
  V(G)$, $c_o(G-S)\leqslant |S|$.
\end{theorem}

\begin{theorem}\label{near_pm}
  A graph $G$ of odd order has no near perfect matching if and only if there exits a set
  $S\subseteq  V(G)$, $c_o(G-S)\geqslant |S|+3$.
\end{theorem}
\begin{proof}
Let $G'$ be a graph obtained from $G$ by adding a new vertex $u$ and
joining $u$ to every vertex of $G$. Then $G$ has a near perfect
matching if and only if $G'$ has a perfect matching.

By Tutte's $1$-Factor Theroem, and the parity, $G'$ has no perfect
matching if and only if there exists a vertex set $S'\subseteq
V(G')$ such that $c_o(G'-S')\geqslant |S'|+2$. Since $u$ is adjacent
to every vertex of $G$, then $u\in S'$. Let $S=S'\setminus \{u\}
\subseteq  V(G)$. Then $c_o(G-S)=c_o(G'-S')\geqslant |S'|+2=|S|+3$.
\end{proof}

The following lemma is proven by Ananchuen and Plummer in
\cite{Ananchuen2007}, they are useful to deal with the graphs with
smaller cut sets. We will use them in our proof several times.

\begin{lemma}\label{cut}
  Let $G$ be a $3$-$\gamma$-vertex-critical graph.
\begin{enumerate}
  \item  If $G$ is disconnected, then $G=3K_{1}$ or $G$ is a disjoint union of a
$2$-$\gamma$-vertex-critical graph and an isolated vertex;
  \item  If $G$ has a cut-vertex $u$, then $c(G-u)=2$. Furthermore, let $C_{i}$ be a
  component of $G-u$ $(i =1, 2)$, then $G[V(C_{i})\cup \{u\}]$ is
  $2$-$\gamma$-vertex-critical;
  \item If $G$ has a $2$-cut $S$, then $c(G-S)\leqslant 3$. Furthermore, if
  $c(G-S)=3$, then $G-S$ must contain at least one singleton.
\end{enumerate}
\end{lemma}

We also need the following results in our proof.

\begin{lemma}[Wang and Yu, \cite{Wang2008}]\label{degree1}
  Let $G$ be a $3$-$\gamma$-vertex-critical graph and $S\subseteq V(G)$.
  If $D_u\subseteq S$ for each vertex $u\in S$, then there exists no
  vertex of degree one in $G[S]$.
\end{lemma}

\begin{theorem}[Wang and Yu, \cite{Wang}]\label{3connected}
  Let $G$ be a $3$-$\gamma$-vertex-critical graph of even order. If
  the minimum degree is at least three, then $G$ is $3$-connected.
\end{theorem}

\begin{theorem}[Mantel, see \cite{West}]\label{triangle-free}
  The maximum number of edges in a triangle-free simple
  graph of order $n$ is $\lfloor {n^2 \over 4}\rfloor$.
\end{theorem}

\section{Proof of Theorem \ref{k167_theorem}}

In this section, we provide a proof of Theorem \ref{k167_theorem}.
\begin{proof}
Suppose, to the contrary, that the theorem does not hold. From
Theorem \ref{tutte_thm} and Theorem \ref{near_pm}, and the parity, there
exists a vertex set $S\subseteq  V(G)$, such that $c_o(G-S)\geqslant
|S|+k-4$ ($k = 6, 7$). Without loss of generality, let $S$ be
\emph{minimal} such a set. By Lemma \ref{cut}, $|S|\geqslant 3$.

\begin{claim}\label{3adjacent}
Each vertex of $S$ is adjacent to at least three odd components of
 $G-S$.
\end{claim}

    Otherwise, there exists a vertex $v\in S$ such that $v$ is adjacent
to at most two odd components of $G-S$. Let $S'=S-\{v\}$. It is easy
to see that $S'$ is a nonempty set which satisfies the condition
$c_o(G-S')\geqslant |S'|+k-4$, contradicting the minimality of
$S$.\\

    Let $C_1,C_2,\ldots,C_t$ be the odd components and
$E_1,E_2,\ldots,E_n$ be the even components of $G-S$.

\vspace{3mm}

\noindent \textbf{Case 1.} $|S|=3$, say $S=\{u,v,w\}$.

Then $t\geqslant k-1$.
\begin{claim}\label{3subset}
  For every vertex $s\in S$, $D_s\subseteq S$.
\end{claim}

Clearly, $D_s\cap S\neq \emptyset$. Assume $D_v=\{u,v'\}$, where
$v'\in V(C_1\cup E_1)$. This means that, if the vertex $v'$ is in the
odd component of $G-S$, we assume $v'\in V(C_1)$; if it is in the
even component of $G-S$, we assume $v'\in V(E_1)$. By Fact 2,
$vu\not\in E(G), vv'\not\in E(G)$, and $u$ dominates $C_2 \cup C_3
\cup \dots \cup C_t$. By Claim \ref{3adjacent}, $w$ is adjacent to at
least two of $C_2, C_3, \dots, C_t$. Without loss of generality, let
$wc_i\in E(G)$, for some $c_i\in V(C_i)$, $i=2,3$. By Fact $2$ again,
$D_{c_i}\cap S=\{v\}$, $i=2,3$. Then $vc_i\notin E(G)$. Since
$vv'\not\in E(G)$, then $D_{c_2}\cap V(C_1 \cup E_1)\neq\emptyset$.
But $D_{c_2}$ can not dominate $c_3$, a contradiction. The claim is
proved.

\vspace{3mm}

By Claim \ref{3subset} and Fact $2$, $S$ is an independent set, and for
any vertex $x\notin S$, $|N_S(x)|\geqslant 2$. In fact, $|N_S(x)|=2$.
Since, if $|N_S(x)|=3$, then $D_x \cap S=\emptyset$.
\begin{claim}
If $t\geqslant 5$, then $G-S$ has no even component.
\end{claim}

Suppose, to the contrary, that there exists an  even component $E_1$.
Choose a vertex $x\in V(E_1)$, and consider $D_x$. Assume
$D_x=\{u,u'\}$, where $u\in S$ and $u'$ is in $C_1$ or in an even
component. Then $u$ dominates $C_2 \cup C_3\cup \dots \cup C_t$. By
Claim \ref{3adjacent}, $w$ is adjacent to at least two of $C_2, C_3,
\dots, C_t$. Without loss of generality, let $wc_i\in E(G)$, where
$c_i\in V(C_i)$, $i=2,3$. By Fact $2$, $D_{c_i}\cap S=\{v\}$, thus
$vc_i\not\in E(G)$ for $i=2,3$. Then $D_{c_2}\cap
V(C_3)\neq\emptyset$ and $v$ dominates $C_1 \cup C_4 \cup C_5 \cup
E_1$. Henceforth $D_{c_j}\cap S=\{w\}$ and $wc_j\not\in E(G)$,
$j=4,5$. Consider $D_{c_4}$, since $wc_5\not\in E(G)$, then
$D_{c_4}\cap V(C_5)\neq\emptyset$ and hence $w$ dominates $C_1 \cup
C_2 \cup C_3$ and $E_1$. Since every vertex of $C_1$ is adjacent to
both $w$ and $v$, then $u$ is not adjacent to any vertex of $C_1$,
hence $u'\in V(C_1)$. Since $\{u,u'\}$ dominates $G - \{x\}$, then
$u$ dominates $E_1-\{x\}$. Since $|E_1|\geqslant 2$, then every
vertex of $V(E_1)-\{x\}$ is adjacent to every vertex of $S$, a
contradiction. So $G-S$ has no even component.

\vspace{3mm}

\textbf{Case 1.1.} There exists a (odd) component, say $C_1$, and a
vertex $c\in V(C_1)$ such that $D_c\cap V(C_1)\neq\emptyset$.

Let $D_c=\{u,c'\}$, where $c'\in V(C_1)$. Then $u$ dominates $C_2
\cup C_3\cup \dots \cup C_t$. Let $c_i \in V(C_i)$, $i=2,\dots, t$.
Since $|N_S(c_i)|=2$ and $uc_i \in E(G)$, assume $wc_2\in E(G)$ and
$wc_3\in E(G)$. Then $D_{c_i}\cap S=\{v\}$ and $vc_i\not\in E(G)$ for
$i=2,3$.  Since $vc_3\not\in E(G)$, then $D_{c_2}\cap
V(C_3)\neq\emptyset$. Therefore, $v$ dominates $C_1 \cup C_4 \cup
C_5$, and hence $wc_4\not\in E(G)$ and $wc_5\not\in E(G)$. Then $w$
dominates $C_1 \cup C_2 \cup C_3$. So every vertex of $C_1$ is
adjacent to both $w$ and $v$, then $u$ is not adjacent to any vertex
of $C_1$. Therefore, for any vertex $x \in V(C_1)$, $D_x \cap S =
\{u\}$ and $|D_x \cap V(C_1)| = 1$. It is easy to see that $C_1$ is
$2$-$\gamma$-vertex-critical, and thus $|V(C_1)|$ is even, a
contradiction.

\vspace{3mm}

\textbf{Case 1.2.} For any vertex $x$ of $C_i$, $D_x\cap
V(C_i)=\emptyset$.

Assume that $|V(C_1)|\geqslant 3$. Let $x\in V(C_1)$, $D_x=\{u,x'\}$.
By Claim $2$ and the assumption $D_x \cap V(C_1) = \emptyset$, we may
assume that $x'\in V(C_2)$. Then $u$ dominates $C_3\cup C_4\cup C_5$
and $C_1-\{x\}$. Since $|N_S(c_i)|=2$ and $uc_i \in E(G)$ for
$i=3,4,5$, so we assume $wc_3 \in E(G)$ and $wc_4 \in E(G)$. Then
$vc_3\not\in E(G)$ and $vc_4\not\in E(G)$. So $D_{c_3}\cap
V(C_4)\neq\emptyset$. It yields that $v$ dominates $C_1$. Since every
vertex of $V(C_1)-\{x\}$ is adjacent to both $u$ and $v$, then it is
not adjacent to $w$. Let $y\in V(C_1)-\{x\}$. Then $D_y \cap
S=\{w\}$, by the assumption $D_y\cap V(C_1)=\emptyset$, so $D_y$ can
not dominate $V(C_1)-\{x,y\}$, a contradiction.

Therefore all the components of $G-S$ are singletons, i.e.,
$C_i=\{c_i\}$. Assume $D_{c_1}=\{u,c_2\}$. Then $uc_1\not\in E(G)$,
$c_2v\in E(G)$ and $c_2w\in E(G)$. Since $|N_S(c_2)|=2$, then
$c_2u\not\in E(G)$. Thus $u$ dominates $G-S-\{c_1, c_2\}$. Therefore,
$D_{c_2}=\{u,c_1\}$. Similarly, we see $D_{c_3}=\{v,c_4\}$,
$D_{c_4}=\{v,c_3\}$, $D_{c_5}=\{w,c_6\}$ and $D_{c_6}=\{w,c_5\}$.
Hence, there is only one 9-vertex graph
satisfying these conditions (see Fig. \ref{9vertex}).\\

\begin{figure}
  \begin{center}
  \includegraphics[width=5cm]{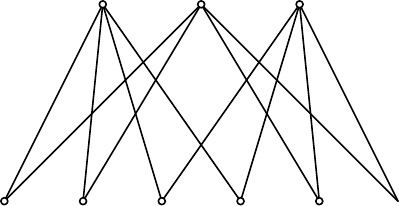}\\
  \caption{A 9-vertex graph which has no near perfect matching.}\label{9vertex}
  \end{center}
\end{figure}

\noindent \textbf{Case 2.} $|S|=4$, and thus $t\geqslant k$.

    We first show that there exists a vertex $a\in S$ such that
$D_a\nsubseteq S$. Otherwise, $D_b\subseteq S$ for every vertex
$b\in S$. By  Fact $2$ and Lemma \ref{degree1}, $S$ is an independent
set. It is easy to check that this is impossible.

So let $u$ be a vertex of $S$ with $D_u\nsubseteq S$. Clearly,
$D_u\cap S\neq \emptyset$. Let $D_u=\{v,x\}$, where $v\in S$ and
$x\in V(G)-S$. Since $G$ is $K_{1,k}$-free, so $t=k$ and $G-S$ has no
even component. Without loss of generality, let $x\in V(C_1)$, then
$v$ dominates all vertices of $\bigcup _{i=2}^{k}V(C_i)$.  Moreover,
 by $K_{1,k}$-freedom again,
$C_2, C_3, \dots, C_k$ are all complete, and $v$ is not adjacent to
any vertex of $V(C_1)$.

Let $S-\{u,v\}=\{w,z\}$. By Claim \ref{3adjacent}, let $wc_i\in E(G)$,
for some $c_i\in V(C_i), i=2,3$. Then $z\in D_{c_2}$. Otherwise, we
have $D_{c_2}\cap S=\{u\}$. Since $ux\notin E(G)$, then $D_{c_2}\cap
V(C_1)\neq \emptyset$, but then $D_{c_2}$ can not dominate $v$, a
contradiction. Similarly, $z\in D_{c_3}$, thus $zc_2\not\in E(G)$
and $zc_3\not\in E(G)$. By Facts 2 and 3, either $D_{c_2}\neq
\{u,z\}$ or $D_{c_3}\neq \{u,z\}$. Assume that $D_{c_2}\neq
\{u,z\}$, thus $D_{c_2}\cap S=\{z\}$. Since $zc_3\not\in E(G)$, then
$D_{c_2}\cap V(C_3)\neq \emptyset$, and $z$ dominates $V(C_1)\cup
V(C_4)\cup V(C_5)\cup V(C_6)$. By a similar argument, $w\in
D_{c_j}$, for some $c_j\in V(C_j), j=4,5,6$. Furthermore,
$wc_j\not\in E(G), j=4,5,6$. From Fact 3, $D_{c_4}\neq\{u,w\}$ or
$D_{c_5}\neq\{u,w\}$ or $D_{c_6}\neq\{u,w\}$. Assume
$D_{c_4}\neq\{u,w\}$. Since $wc_5\not\in E(G)$, then $D_{c_4}\cap
V(C_5)\neq\emptyset$, but $D_{c_4}$ can not dominate $c_6$, a
contradiction.

\vspace{3mm}

\noindent \textbf{Case 3.} $|S|=5$, and thus $t\geqslant k+1$.

\begin{claim}\label{5subset}
For every vertex $s\in S$, $D_s\subseteq S$.
\end{claim}

Otherwise, $D_u\nsubseteq S$ for some $u\in S$. Clearly, $D_u\cap
S\neq \emptyset$. Let $D_u=\{y,z\}$, where $y\in S$ and $z\not\in S$.
Since $t\geqslant k+1$, $y$ must dominate at least $k$ odd components
of $G-S$, which contradicts to $K_{1,k}$-freedom.

By Claim \ref{5subset} and Lemma \ref{degree1}, each vertex of $S$  has
degree $0$ or $2$ in $G[S]$. It is not hard to see that $G[S]$ can
only be a $5$-cycle or a disjoint union of a $4$-cycle and an
isolated vertex. Let $S=\{s_1,s_2,s_3,s_4,s_5\}$. There are
$\binom{5}{2}=10$ distinct pairs of vertices in $S$. By Fact $3$ and
Claim \ref{5subset}, there must exist a vertex $x$ in an odd component
of $G-S$ such that $D_x\nsubseteq S$. Assume that $x\in V(C_1)$.
Clearly, $D_x\cap S\neq \emptyset$. Since $G$ is $K_{1,k}$-free, we
have $t=k+1$ and $G-S$ has no even component.

\vspace{3mm}

\textbf{Case 3.1.} $G[S]$ is a $5$-cycle.

Let $s_1s_2s_3s_4s_5s_1$ be the $5$-cycle in the counterclockwise
order and $D_x=\{s_1,x'\}$, where $x'\not\in S$. Since $G$ is
$K_{1,k}$-free, then $x'\notin V(C_1)$. Assume that $x'\in V(C_2)$.
Then $s_1$ dominates $\bigcup _{i=3}^{k+1} V(C_i)$ and $x'$ is
adjacent to both $s_3$ and $s_4$. Moreover,  $K_{1,k}$-freedom of $G$
implies that $C_3, C_4, \dots, C_{k+1}$ are all complete and $s_1$ is
not adjacent to any vertex of $V(C_1)\cup V(C_2)$. Henceforth, $C_1$
is a singleton (i.e., $V(C_1)=\{x\}$).

Since $D_{s_3}=\{s_1,s_5\}$, then $s_5$ dominates $V(C_1)\cup
V(C_2)$. Similarly, since $D_{s_4}=\{s_1,s_2\}$, $s_2$ dominates
$V(C_1)\cup V(C_2)$. Therefore, $x'$ is adjacent to all vertices of
$S-\{s_1\}$. Now consider $D_{x'}$. Since $D_{x'}\cap S = \{s_1\}$
and $s_1x\not\in E(G)$, it follows that $D_{x'}=\{s_1,x\}$. Hence,
$x$ is adjacent to both $s_3$ and $s_4$, and $V(C_2)=\{x'\}$. But
then $\{s_1,s_3\}$ is a dominating set in $G$, a contradiction to
$\gamma(G)=3$.

\vspace{3mm}

\textbf{Case 3.2.} $G[S]$ is a disjoint union of a $4$-cycle and an
isolated vertex.

Let $s_1s_2s_3s_4s_1$ be the $4$-cycle in the counterclockwise order
and $s_5$ be the isolated vertex in $G[S]$. Then
$D_{s_1}=\{s_3,s_5\}$, $D_{s_2}=\{s_4,s_5\}$, $D_{s_3}=\{s_1,s_5\}$,
and $D_{s_4}=\{s_2,s_5\}$.

Since $G$ is $K_{1,k}$-free, $s_5$ is adjacent to at most $k-1$ (odd)
components of $G-S$. Without loss of generality, let $C_1,\ldots,C_r$
be the components which are not adjacent to $s_5$. Then $t=k+1$
implies $r\geqslant 2$. Thus $s_i$ dominates $\bigcup _{j=1}^rV(C_j)$
for $i=1,2,3,4$. Now consider $D_{c_1}$, where $c_1 \in V(C_1)$.
Clearly, $D_{c_1} \cap S = \{s_5\}$. Since $s_5$ is not adjacent to
$V(C_2)$, then $D_{c_1}\cap V(C_2)\neq\emptyset$. Therefore, $r=2$
and $s_5$ dominates $\bigcup _{j=3}^{k+1} V(C_j)$. Moreover,
$V(C_1)=\{c_1\}$. By a similar argument, $C_2$ is also a singleton.

    For any vertex $v\in\bigcup _{j=3}^{k+1} V(C_j)$, by Fact 2,
$s_5\not\in D_{v}$, but the vertices in $S-\{s_5\}$ do not dominate
$s_5$. Then $D_v\not\subseteq S$ and $D_v\cap
\{s_1,s_2,s_3,s_4\}\neq\emptyset$. From $K_{1,k}$-freedom of $G$, it
implies that $C_3, C_4, \dots, C_{k+1}$ are all singletons, say
$V(C_j)=\{c_j\}$ for $j=3,\dots, k+1$. Then $|V(G)|=12$ or $13$ (see
examples: Fig. \ref{k16}, Fig. \ref{k17}).\\

\begin{figure}[h]
\begin{minipage}[t]{0.45\linewidth}
\centering
\includegraphics[width=\textwidth]{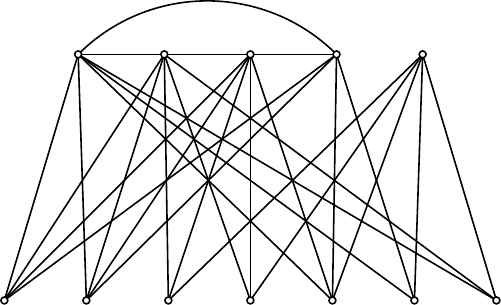}
\caption{A $K_{1,6}$-free graph without perfect matching.
\label{k16}}
\end{minipage}
\hfill
\begin{minipage}[t]{0.45\linewidth}
\centering
\includegraphics[width=\textwidth]{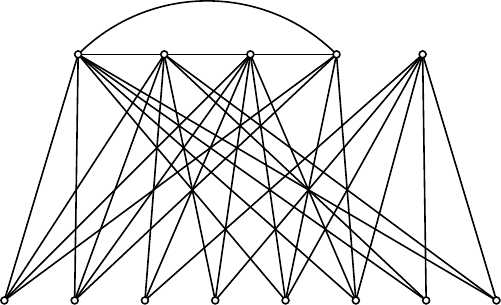}
\caption{A $K_{1,7}$-free graph without near perfect
matching.\label{k17}}
\end{minipage}
\end{figure}

\noindent \textbf{Case 4.} $|S|\geqslant 6$, and thus $t\geqslant
k+2$.

\begin{claim}\label{6subset}
 For every vertex $s\in V(G)$, $D_s\subseteq S$.
\end{claim}
    Suppose that $D_x\nsubseteq S$ for some $x\in V(G)$. Clearly,
$D_x\cap S\neq \emptyset$. Let $D_x=\{y,z\}$, where $y\in S$ and
$z\not\in S$. Since $t\geqslant k+2$, $y$ must dominate at least $k$
odd components of $G-S$, a contradiction.

For each $i=1, \dots, t$, let $S_i\subseteq S$ be the set of vertices
in $S$ which are adjacent to some vertex in $C_i$, and let
$d=\min\{|S_i|\}$. Without loss of generality, assume that $|S_1|=d$.
Note that for any vertex $v\in V(G)-V(C_1)$, $D_v \subset S$ has to
dominate $C_1$, thus, $D_v\cap S_1\neq \emptyset$. We call such a set
$D_v$ \emph{normal $2$-set associated with $v$ and $S_1$}, or
\emph{normal set} in short. By a simple counting argument, we see
that there are at most $\binom{|S|}{2}-\binom{|S|-d}{2}$ normal sets.

\vspace{3mm}

\textbf{Case 4.1.} $G$ is $K_{1,6}$-free, and $|V(G)|$ is even.

    Since every vertex in $S$ is adjacent to at most five components of
$G-S$, then $c(G-S)\leqslant 10$. Henceforth, $6\leqslant
|S|\leqslant 8$ and $d\leqslant \lfloor\frac{5|S|}{|S|+2}\rfloor
\leqslant 4$.

If $|S|=6$, then $\binom{6}{2}-\binom{6-d}{2}\geqslant 13$, and thus
$d\geqslant 4$. But $d\leqslant \lfloor\frac{5 \times 6}{6+2}\rfloor
<4$, a contradiction.

If $|S|=7$, then
    \begin{equation}\label{eq1}
      \binom{7}{2}-\binom{7-d}{2}\geqslant 15
    \end{equation} or
$d\geqslant 3$. Since $d\leqslant \lfloor \frac{5 \times
7}{7+2}\rfloor <4$, then $d=3$ and the equality holds in (\ref{eq1}).
Let $S_1=\{u,v,w\}$, then $\{u,v\}$, $\{u,w\}$, $\{v,w\}$ are all
corresponding to some $D_x$ where $x\not\in V(C_1)$. Since $u$ is
adjacent to at most five components of $G-S$, so we may assume that
$u$ is not adjacent to $C_6, C_7, \dots, C_9$. Then $v$ dominates at
least three of them, and $v$ is adjacent to at most two of $C_1, C_2,
\dots, C_5$. Similarly, $w$ is adjacent to at most two of $C_1, C_2,
\dots, C_5$. Both $v$ and $w$ are adjacent to $C_1$, then $\{v,w\}$
can dominate at most two of $C_2, C_3, \dots, C_5$, hence it can not
be realized a $D_x$ for some $x \not\in V(C_1)$, a contradiction.

If $|S|=8$, then $c(G-S)=c_o(G-S)=10$. We construct a graph $H$ with
vertex set $S$ and $uv\in E(H)$ if and only if $D_x=\{u,v\}$ for some
$x\in V(G)$. We show that $H$ is triangle-free. Let $u,v,w\in S$, if
$uv\in E(H)$, $uw\in E(H)$ and $u$ is not adjacent to $C_6, \dots,
C_{10}$, then both $v$ and $w$ are adjacent to at least four of them.
Hence both $v$ and $w$ are adjacent to at most one component of $C_1,
C_2, \dots,C_5$. Therefore $\{v,w\}$ is not a $D_x$ for any $x\in
V(G)$. By Theorem \ref{triangle-free}, $|E(H)|\leqslant
\lfloor\frac{8^2}{4}\rfloor=16< |V(G)|$, a contradiction.

\vspace{3mm}

\textbf{Case 4.2.} $G$ is $K_{1,7}$-free and $|V(G)|$ is odd.

Since every vertex in $S$ is adjacent to at most six components of
$G-S$, then $c(G-S)\leqslant 12$. So $6 \leqslant |S|\leqslant 9$.

If $|S|=6$, by Claim \ref{6subset} and Fact $3$, $\binom{6}{2}\geqslant
|V(G)| \geqslant 6 + 9$. Then $|V(G)|=15$, and $G-S$ is an
independent set of nine vertices. Moreover, every pair in $S$ is
corresponding to a $D_x$ for some $x\in V(G)$. As
$\binom{6}{2}-\binom{6-d}{2}\geqslant 14$, so $d\geqslant 4$. For any
$x\not\in S$, $D_x \subset S$, by Fact 2, every vertex in $G-S$ has
degree 4, and then every vertex of $S$ is adjacent to six components
of $G-S$. Let $\delta$ be the minimum degree of $G[S]$ and 
$d_{G[S]}(u) = d$. If $d \leq 2$, then there exists at least one pair
in $S \setminus N_{G[S]}[u]$ which is not corresponding to $D_u$, and
thus it does not dominate $u$, a contradiction. By Fact 2, $G[S]$ is
a $3$-regular graph. From the above information, it is not hard to
see that there are only two such graphs (see Fig. \ref{s45}).

\begin{figure}[h]
\begin{minipage}[t]{0.45\linewidth}
\centering
\includegraphics[width=\textwidth]{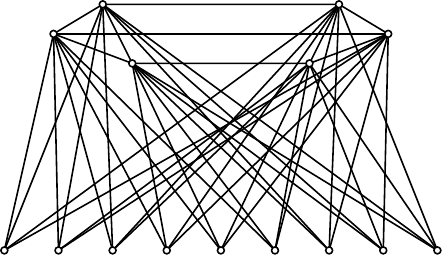}
\end{minipage}
\hfill
\begin{minipage}[t]{0.45\linewidth}
\centering
\includegraphics[width=\textwidth]{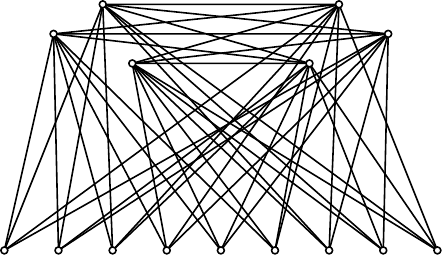}
\end{minipage}
\caption{Two exceptions when $G$ is $K_{1,7}$-free and $|V(G)|$ is
odd. \label{s45}}
\end{figure}

If $|S|=7$, we construct an auxiliary graph $H$ with vertex set $S$
and $uv\in E(H)$ if and only if $D_x=\{u,v\}$ for some $x\in S$.
Assume that $uv, uw \in E(H)$, and $u$ is not adjacent to $C_7,
\dots, C_{10}$. Then both $v$ and $w$ dominate $C_7, \dots, C_{10}$,
and are all adjacent to at most two of $C_1, C_2, \dots, C_6$. Hence
$\{v,w\}$ can not be realized as a $D_v$ for some $v\in V(G)$.
Therefore, $H$ is triangle-free. If $H$ contains a cycle of length
at least five, then at least five pairs can not be realized as a
$D_x$ for some $x\in V(G)$, $\binom{7}{2}-5=16< 17\leqslant |V(G)|$,
a contradiction. As $|E(H)|>|V(H)|-1$, so $H$ only contains cycles
of length four, and $H$ is bipartite. Let $s_1s_2s_3s_4$ be a four
cycle in $H$. $|E(H)|>|V(H)|-1=6$, it yields that the component
which contains the $4$-cycle $s_1s_2s_3s_4$, say $H'$, has at least
six vertices. The pairs in the same partite of $H'$ can not be
realized as a $D_x$ for some $x\in V(G)$, a simple counting argument
shows that $H$ has at least five such pairs. So $\binom{7}{2}-5=16<
17\leqslant |V(G)|$, a contradiction.

If $8\leqslant |S|\leqslant 9$, we construct a graph $H$ as in the
case that ``$G$ is $K_{1,6}$-free, $|V(G)|$ is even, and $|S|=8$''.
Similarly, $H$ is triangle-free, by Theorem \ref{triangle-free},
$|E(G)|\leqslant \lfloor\frac{|S|^2}{4}\rfloor<|V(G)|$, a
contradiction.
\end{proof}

\vspace{4mm}

\begin{remark} The conclusion in this theorem holds for all graphs
except $|V(G)|=12$ or $13$. For these cases, we can determine the
exceptions precisely in some cases (such as in Case 4.2) but fail to
determine all of them in other cases (such as in Case 3.2). With
some efforts, one may be able to find all graphs which have no
perfect matching or near-perfect matching for $|V(G)|=12$ or $13$.
\end{remark}

\begin{remark}
Ananchuen and Plummer \cite{Ananchuen2006} showed that: let $G$ be a
connected $3$-$\gamma$-vertex-critical graph of even order. If G is
claw-free, then G is bicritical. The authors also generalized this
result, and proved that: let $G$ be a $3$-$\gamma$-vertex-critical
graph of even order, if $G$ is $K_{1,4}$-free, and the minimum
degree is at least four, then $G$ is bicritical. This result will be
published in a future article.
\end{remark}

{\bf Acknowledgment } The authors would like to thank the anonymous
referee and the editor for their valuable comments and suggestions.


\begin{thebibliography}{1}
\expandafter\ifx\csname url\endcsname\relax
  \def\url#1{\texttt{#1}}\fi
\expandafter\ifx\csname urlprefix\endcsname\relax\def\urlprefix{URL
}\fi \expandafter\ifx\csname href\endcsname\relax
  \def\href#1#2{#2} \def\path#1{#1}\fi

\bibitem{Ananchuen2007}
N.~Ananchuen, M.~D. Plummer, {Matchings in 3-vertex-critical graphs:
The odd case}, Discrete Mathematics 307~(13) (2007) 1651--1658.

\bibitem{Ananchuen2006}
N.~Ananchuen, M.~D. Plummer, On the connectivity and matchings in
  3-vertex-critical claw-free graphs, Util. Math. 69 (2006) 85--96.

\bibitem{Ananchuen2005}
N.~Ananchuen, M.~D. Plummer, {Matchings in 3-vertex-critical
  graphs: The even case}, Networks 45~(4) (2005) 210--213.

\bibitem{Brigham1988}
R.~C. Brigham, P.~Z. Chinn, R.~D. Dutton, {Vertex
domination-critical
  graphs}, Networks 18~(3) (1988) 173--179.

\bibitem{Flandrin1999}
E.~Flandrin, F.~Tian, B.~Wei, L.~Zhang, {Some properties of
3-domination-critical graphs},
  Discrete Mathematics 205~(1-3) (1999) 65--76.

\bibitem{Fulman1995}
J.~Fulman, D.~Hanson, G.~MacGillivray, {Vertex domination-critical
  graphs}, Networks 25~(2) (1995) 41--43.

\bibitem{Lov'asz1986}
L.~Lov\'{a}sz, M.~D. Plummer, Matching Theory, North-Holland Inc.,
Amsterdam, 1986.

\bibitem{Wang2008}
T. Wang, Q. Yu, Factor-critical property in 3-dominating-critical
graphs, Discrete Math. (2008), doi:10.1016/j.disc.2007.11.062

\bibitem{Wang}
T.~Wang, Q.~Yu, A conjecture on $k$-factor-critical and
$3$-$\gamma$-critical graphs (submitted).

\bibitem{West}
D. B.~West, Introduction to Graph Theory, 2nd Edition, Prentice
Hall, 2000.


\end{thebibliography}
\end{document}